\documentclass{amsart}

\usepackage{amssymb}
\usepackage{amsmath}
\usepackage{mathrsfs}
\usepackage{dsfont}
\usepackage{url}

\newtheorem*{theorem*}{Theorem}
\newtheorem*{lemma*}{Lemma}

\newtheorem*{prime-k-tuple-conjecture}{Conjecture H}

\newtheorem*{acknowledgment}{Acknowledgment}

\begin{document}

\title[On the differences between consecutive prime numbers, I]{On the differences between consecutive prime numbers, I}

\author[D. A. Goldston]{D. A. Goldston$^{*}$}

\address{Department of Mathematics,
 San Jos\'{e} State University,
 315 MacQuarrie Hall,
 One Washington Square,
 San Jos\'{e}, California 95192-0103,
 United States of America}

\email{daniel.goldston@sjsu.edu}

\thanks{$^{*}$During the preparation of this work, the first author received support from the National Science Foundation Grant DMS-1104434.}

\author[A. H. Ledoan]{A. H. Ledoan}

\address{Department of Mathematics,
 University of Tennessee at Chattanooga,
 417F EMCS Building (Department 6956),
 615 McCallie Avenue,
 Chattanooga, Tennessee 37403-2598,
 United States of America}

\email{andrew-ledoan@utc.edu}



\begin{abstract}
We show by an inclusion-exclusion argument that the prime $k$-tuple conjecture of Hardy and Littlewood provides an asymptotic formula for the number of consecutive prime numbers which are a specified distance apart. This refines one aspect of a theorem of Gallagher that the prime $k$-tuple conjecture implies that the prime numbers are distributed in a Poisson distribution around their average spacing.
\end{abstract}

\maketitle

\thispagestyle{empty}

\section{Introduction and statement of results}

In 1976, Gallagher \cite{Gallagher1976}, \cite{Gallagher1981} showed that a uniform version of the prime $k$-tuple conjecture of Hardy and Littlewood implies that the prime numbers are distributed in a Poisson distribution around their average spacing. Specifically, let $P_{r} (h, N)$ denote the number of positive integers $n$ less than or equal to $N$ such that the interval $(n, n + h]$ contains exactly $r$ prime numbers. Gallagher then proved that an appropriate form of the prime $k$-tuple conjecture implies, for  any positive constant $\lambda$ and $h \sim \lambda \log N$ as $N \to \infty$, that
\[
P_{r} (h, N)
 \sim e^{-\lambda} \frac{\lambda^{r}}{r!} N.
\]

In particular, if $r = 0$, then we obtain by an argument using the prime number theorem that, as $N \to \infty$,
\begin{equation}\label{eq1}
\sum_{\substack{p_{n + 1} \leq N \\ p_{n+1} - p_{n} \geq \lambda \log n}} 1
 \sim e^{-\lambda} \frac{N}{\log N}.
\end{equation}
Here, $p_{n}$ is used to denote the $n$th prime number. The purpose of the present paper is to obtain a refinement of \eqref{eq1}, which shows that the Poisson distribution of the prime numbers in  short intervals extends down to the individual differences between consecutive prime numbers. To obtain this result, we employ a version of the prime $k$-tuple conjecture formulated as Conjecture H in Section \ref{section2}, which is equivalent to the form of the conjecture used by Gallagher.

\begin{theorem*}
Assume Conjecture H. Let $d$ be any positive integer, and let $p$ be a prime number. Let, further,
\[
\mathfrak{S}(d)
 = \left\{ \begin{array}{ll}
      {\displaystyle 2 C_{2} \prod_{\substack{ p \mid d \\ p > 2}} \left(\frac{p - 1}{p - 2}\right),} & \mbox{if $d$ is  even;} \\
      0, & \mbox{if $d$ is  odd;} \\
\end{array}
\right.
\]
where
\[
C_{2}
 = \prod_{p > 2} \left( 1 - \frac{1}{(p - 1)^{2}}\right)
 = 0.66016\ldots,
\]
and define
\[
N(x, d)
 = \sum_{\substack{ p_{n + 1} \leq x \\ p_{n + 1} - p_{n} = d}} 1,
\]
where $p_{n}$ denotes the $n$th prime number. Then for any positive constant $\lambda$ and $d$ even with $d \sim \lambda \log x$ as $x \to \infty$, we have
\begin{equation}\label{eq2}
N(x, d)
 \sim e^{-\lambda} \mathfrak{S}(d) \frac{x}{(\log x)^2}.
\end{equation}
\end{theorem*}

Here, we note that $\mathfrak{S}(d)$ is the singular series in the conjectured asymptotic formula for the number of prime pairs differing by $d$. Our theorem shows that, for consecutive prime numbers, the Poisson density  is superimposed onto this formula for prime pairs.

Our theorem as well as its proof are implicitly contained in the 1999 paper of Odlyzko, Rubinstein and Wolf \cite{OdlyzkoRubinsteinWolf1999} on jumping champions.\footnote{An integer $d$ is called a jumping champion for a given $x$ if $d$ is the most frequently occurring difference between consecutive prime numbers up to $x$.} Without claiming any originality, we think it is worthwhile to explicitly state and prove \eqref{eq2}. More precise results when $d / \log x \to 0$ will be addressed in a second paper.

\section{The Hardy-Littlewood prime $k$-tuple conjectures}\label{section2}

Let $\mathcal{H} = \{h_{1}, \ldots , h_{k}\}$ be a set of $k$ distinct integers. Let $\pi(x; \mathcal{H})$ denote the number of positive integers $n$ less than or equal to $x$ for which $n + h_{1}, \ldots , n + h_{k}$ are simultaneously prime numbers. Then the prime $k$-tuple conjecture of Hardy and Littlewood \cite{HardyLittlewood1922} is that, for $x \to \infty$,
\begin{equation}\label{eq3}
\pi(x; \mathcal{H})
 \sim \mathfrak{S}(\mathcal{H}) \ \! \mbox{\textup{li}}_{k}(x),
\end{equation}
where
\[
\mathfrak{S}(\mathcal{H})
 = \prod_{p} \left(1 - \frac{1}{p}\right)^{-k} \left(1 - \frac{\nu_{\mathcal{H}}(p)}{p}\right),
\]
$\nu_{\mathcal{H}}(p)$ denotes the number of distinct residue classes modulo $p$ occupied by the elements of $\mathcal{H}$, and
\begin{equation}\label{eq4}
\mbox{\textup{li}}_{k}(x)
 = \int_{2}^{x} \frac{\,dt}{(\log t)^{k}}.
\end{equation}

Note in particular that, if $\nu_{\mathcal{H}}(p) = p$ for some prime number $p$, then   $\mathfrak{S}(\mathcal{H}) = 0$. However, if $\nu_{\mathcal{H}}(p) < p$ for all prime numbers $p$, then  $\mathfrak{S}(\mathcal{H})\neq 0$ in which case the set $\mathcal{H}$ is called {\it admissible}. In \eqref{eq3}, $\mathcal{H}$ is assumed to be admissible, since otherwise $\pi(x; \mathcal{H})$ is equal to 0 or 1.

The prime $k$-tuple conjecture has been verified only for the prime number theorem. That is to say, for the case of $k = 1$. It has been asserted that, in its strongest form, the conjecture holds true for any fixed integer $k$ with an error term that is $O_{k}(x^{1 / 2 + \varepsilon})$ at most and uniformly for $\mathcal{H} \subset [1, x]$. (See Montgomery and Soundararajan \cite{MontgomerySoundararajan1999-2002}, \cite{MontgomerySoundararajan2004}.) However, we do not need such strong conjectures here. Using
\begin{equation}\label{eq5}
\mbox{\textup{li}}_{k}(x)
 = \frac{x}{(\log x)^{k}} + O\left(\frac{kx}{(\log x)^{k + 1}}\right),
\end{equation}
obtained from integration by parts, we replace $\mbox{\textup{li}}_{k}(x)$ by its main term and make the following conjecture.
\begin{prime-k-tuple-conjecture}
For each fixed integer $k \geq 2$ and admissible set $\mathcal{H}$, we have
\[
\pi(x; \mathcal{H})
 = \mathfrak{S}(\mathcal{H}) \frac{x}{(\log x)^{k}} (1 + o_{k}(1)),
\]
uniformly for $\mathcal{H} \subset [1, h]$, where $h \sim \lambda \log x$ as $x \to \infty$ and $\lambda$ is a positive constant.
\end{prime-k-tuple-conjecture}

\section{Inclusion-exclusion for consecutive prime numbers}

The prime $k$-tuple conjecture for the case when $k = 2$ provides an asymptotic formula for the number of prime numbers with a given difference $d$. We need to find a corresponding formula where we restrict the count to prime numbers that are consecutive, and for this one can use the prime $k$-tuple conjecture with $k = 3, 4, \ldots$ and inclusion-exclusion to obtain upper and lower bounds for the number of consecutive prime numbers with difference $d$. This method has appeared in a series of papers of Brent \cite{Brent1974}, \cite{Brent1975}, \cite{Brent1980} and was used by  Erd\H{o}s and Strauss \cite{ErdosStraus1980} and Odlyzko, Rubinstein and Wolf \cite{OdlyzkoRubinsteinWolf1999} in their study of jumping champions.

We consider a special type of tuple $\mathcal{D}_{k}$ for which
\[
\mathcal{D}_{2}
 = \{0, d\}
\]
and, for $k \geq 3$,
\[
\mathcal{D}_{k}
 = \{0, d_{1}, \ldots , d_{k - 2}, d\}.
\]
Here, we require that $d$ is even. We want to count the number of consecutive prime numbers which do not exceed $x$ and have difference $d$, namely $N(x,d)$, and for this we do inclusion-exclusion with
\[
\pi_{2}(x, d)
 =\sum_{\substack{p \leq x \\ p - p' = d}} 1,
\]
where $p'$ is also a prime number and, for $k \geq 3$,
\[
\pi_{k}(x, d_{1}, \ldots, d_{k - 2}, d)
 = \sum_{\substack{p \leq x \\ p - p' = d \\ p - p_{j} = d_{j}, \!\ 1 \leq j \leq k - 2}} 1.
\]
Inserting the expected main term, we obtain
\begin{equation}\label{eq6}
\pi_{2}(x, d)
 = \mathfrak{S}(d) \ \! \mbox{\textup{li}}_{2}(x) + R_{2}(x, d)
\end{equation}
and, for $k \geq 3$,
\begin{equation}\label{eq7}
\pi_{k}(x, d_{1}, \ldots, d_{k - 2}, d)
 = \mathfrak{S}(\mathcal{D}_{k}) \ \! \mbox{\textup{li}}_{k}(x) + R_{k}(x, \mathcal{D}_{k}).
\end{equation}

We now carry out the inclusion-exclusion. We trivially have
\[
N(x, d)
 \leq \pi_{2}(x, d).
\]
The consecutive prime numbers that differ by $d$ are those prime numbers $p$ and $p'$ satisfying $p - p' = d$ such that there is no third prime number $p''$ with $p' < p'' <p$.  We can exclude these non-consecutive prime numbers differing by $d$  by removing all triples of this form, although this will exclude the same non-consecutive pair of prime numbers more than once if there are quadruples of prime numbers such that $p' < p'' < p''' < p$. Hence, writing $p - p'' = d'$, we obtain the lower bound
\[
N(x, d)
 \geq \pi_{2}(x, d) - \sum_{1 \leq d' < d} \pi_{3}(x, d', d).
\]
We next obtain an upper bound by including the quadruples eliminated in the previous step and continue in this fashion to get, for $R \geq 1$,
\begin{equation}\label{eq8}
Q_{2R + 1}(x, d)
 \leq N(x, d)
 \leq Q_{2R}(x, d),
\end{equation}
where, for $N \geq 2$,
\[
Q_{N}(x, d)
 = \pi_{2}(x, d) + \sum_{k = 3}^{N} (-1)^{k} \sum_{1 \leq d_{1} < \ldots < d_{k - 2} < d} \pi_{k}(x, d_{1}, \ldots , d_{k - 2}, d).
\]
We use the convention here that an empty sum has the value zero.

To evaluate $Q_{N}(x, d)$, we require a special type of singular series average considered by Odlyzko, Rubinstein and Wolf. Let, for $k\geq 3$,
\begin{equation}\label{eq9}
A_{k}(d)
 = \sum_{1 \leq d_{1} < \ldots < d_{k - 2} < d} \mathfrak{S}(\mathcal{D}_{k}).
\end{equation}
Odlyzko, Rubinstein and Wolf \cite{OdlyzkoRubinsteinWolf1999} proved that, for $k \geq 3$,
\begin{equation}\label{eq10}
A_{k}(d)
 = \mathfrak{S}(d) \frac{d^{k - 2}}{(k - 2)!} + E_{k}(d),
\end{equation}
where
\begin{equation}\label{eq11}
E_{k}(d)
 = O_{k} \left(\frac{d^{k - 2}}{\log\log d}\right).
\end{equation}
(See, also, Goldston and Ledoan \cite{GoldstonLedoan2011}.) Thus, on substituting \eqref{eq6}, \eqref{eq7}, \eqref{eq9} and \eqref{eq10}, we find that, for $N \geq 2$,
\[
\begin{split}
Q_{N}(x, d)
 &= \mathfrak{S}(d) \ \! \mbox{\textup{li}}_{2}(x) +  \sum_{k = 3}^{N} (-1)^{k} A_{k}(d) \ \! \mbox{\textup{li}}_{k}(x) + R_{2}(x, d) \\ &\quad +  \sum_{k = 3}^{N} (-1)^{k} \sum_{1 \leq d_{1} < \ldots < d_{k - 2} < d} R_{k}(x, \mathcal{D}_{k}) \\
 &= \mathfrak{S}(d) \int_{2}^{x} \left[\sum_{k = 0}^{N - 2} \frac{1}{k!} \left(\frac{-d}{\log t}\right)^{k}\right] \frac{\,dt}{(\log t)^{2}} + \int_{2}^{x} \sum_{k = 3}^{N} (-1)^{k} E_{k}(d) \frac{\,dt}{(\log t)^{k}} \\ & \quad + R_{2}(x, d) + \sum_{k = 3}^{N} (-1)^{k} \sum_{1 \leq d_{1} < \ldots < d_{k - 2} < d} R_{k}(x, \mathcal{D}_{k}),
\end{split}
\]
where we used \eqref{eq4} in the second line.

We can extract a main term independent of $N$ out of the first term on the far right-hand side above by using Taylor's theorem. With the remainder expressed in Lagrange's form, we have that, for $M \geq 0$ and $x > 0$,
\[
e^{-x}
 = \sum_{k = 0}^{M} \frac{1}{k!} (-x)^{k} + \frac{e^{-\xi}}{(M + 1)!} (-x)^{M + 1},
\]
where $\xi$ lies in the open interval joining 0 and $x$. Hence, we have
\[
\begin{split}
\int_{2}^{x} \left[\sum_{k = 0}^{N - 2} \frac{1}{k!} \left(\frac{-d}{\log t}\right)^{k}\right] \frac{\,dt}{(\log t)^{2}}
 &= \int_{2}^{x} \exp\left(\frac{-d}{\log t}\right) \frac{\,dt}{(\log t)^{2}} \\ &\quad + O\left(\int_{2}^{x} \frac{1}{(N - 1)!} \left(\frac{d}{\log t}\right)^{N - 1} \frac{\,dt}{(\log t)^{2}}\right) \\
 & = \int_{2}^{x} \exp\left(\frac{-d}{\log t}\right) \frac{\,dt}{(\log t)^{2}} \\ &\quad + O\left( \frac{1}{\sqrt{N}} \left(\frac{3d}{N\log x}\right)^{N - 1} \frac{x}{(\log x)^{2}}\right),
\end{split}
\]
by the Stirling formula
\[
(M - 1)!
 = \sqrt{2 \pi} M^{M - 1 / 2} e^{-M} \left[1 + O\left(\frac{1}{M}\right)\right].
\]

Therefore, we have proved the following lemma.

\begin{lemma*}
For $N \geq 2$, we have 
\[
\begin{split}
Q_{N}(x, d)
&= \mathfrak{S}(d) I(x, d) + \int_{2}^{x} \sum_{k = 3}^{N} (-1)^{k} E_{k}(d) \frac{\,dt}{(\log t)^{k}} + R_{2}(x, d) \\ & \quad + \sum_{k = 3}^{N} (-1)^{k} \sum_{1 \leq d_{1} < \ldots < d_{k - 2} < d}R_{k}(x, \mathcal{D}_{k}) \\ &\quad + O\left(\frac{1}{\sqrt{N}} \left(\frac{3d}{N \log x}\right)^{N - 1} \frac{x}{(\log x)^{2}}\right),
\end{split}
\]
where 
\[
I(x, d)
 = \int_{2}^{x} \exp\left(\frac{-d}{\log t}\right) \frac{\,dt}{(\log t)^{2}}.
\]
\end{lemma*}

\section{Proof of the theorem} 

If we had imposed the additional condition that $n + h_{j}$ is less than or equal to $x$, for $j \in \{1, \ldots, k\}$, in the definition of $\pi(x; \mathcal{H})$ in Section \ref{section2}, we would have that $\pi_{2}(n,d)=\pi(x;\mathcal{D}_{2})$ and, for $k \geq 3$, $\pi_{k}(x, d_{1}, \ldots, d_{k - 2}, d) = \pi(x; \mathcal{D}_{k})$. However, this condition has no effect on Conjecture H, since with $\mathcal{H} \subset [1, h]$ the condition removes at most $h$ tuples, which are absorbed into the error term. Thus, assuming that Conjecture H holds true for $k \leq N$, we have that $R_{2}(x, d) = o(\mathfrak{S}(d) x / (\log x)^{2})$, $R_{k}(x, \mathcal{D}_{k}) = o_{k}(\mathfrak{S}(\mathcal{D}_{k}) x / (\log x)^{k})$ and, by \eqref{eq11}, $E_{k}(d) = o_{k}(d^{k - 2})$. Then by the lemma and since $\mathfrak{S}(d) \gg 1$ for $d$ even, we have, for $x \to \infty$ and $d \sim \lambda \log x$,
\[
\begin{split}
Q_{N}(x, d)
 &= \mathfrak{S}(d) I(x, d) + o_{N}\left(e^{2\lambda} \frac{x}{(\log x)^{2}}\right) + o\left(\mathfrak{S}(d) \frac{x}{(\log x)^{2}}\right) \\ & \quad + \sum_{k = 3}^{N} \sum_{1 \leq d_{1} < \ldots < d_{k - 2} < d} o_{k}\left(\mathfrak{S}(\mathcal{D}_{k}) \frac{x}{(\log x)^{k}}\right) \\ &\quad + O\left(\frac{1}{\sqrt{N}} \left(\frac{3\lambda}{N}\right)^{N - 1} \frac{x}{(\log x)^{2}}\right) \\
&= \mathfrak{S}(d) I(x, d) + o_{N}\left(e^{2\lambda} \mathfrak{S}(d) \frac{x}{(\log x)^2}\right) \\ &\quad + O\left(\left(\frac{4\lambda}{N}\right)^{N} \frac{x}{(\log x)^{2}}\right).
\end{split}
\]
Finally, on letting $N$ tend to infinity sufficiently slowly, the theorem follows from the estimate
\[
I(x, d)
 \sim e^{-\lambda} \frac{x}{(\log x)^{2}}
\]
and \eqref{eq8}.

To prove this last estimate, we let  $\bar{d} = d / \log x$ and apply \eqref{eq5} to obtain the upper bound
\[
I(x, d)
 \leq e^{-\bar{d}} \int_{2}^{x} \frac{\,dt}{(\log t)^{2}}
 = e^{-\bar{d}} \left[\frac{x}{(\log x)^{2}} + O\left(\frac{x}{(\log x)^{3}}\right)\right],
\]
and the lower bound
\[
\begin{split}
I(x, d)
 &\geq \exp\left(\frac{-d}{\log x - \log\log x}\right) \int_{x / \log x}^{x} \frac{dt}{(\log t)^{2}} \\
 &\geq \exp\left[{-\bar{d} \left(1 + O\left(\frac{\log\log x}{\log x}\right)\right)}\right] \left[\mbox{\textup{li}}_{2}(x) - \mbox{\textup{li}}_{2} \left(\frac{x}{\log x}\right)\right] \\
 &\geq e^{-\bar{d}} \frac{x}{(\log x)^{2}} \left[1 + O\left(\frac{\bar{d} \log\log x}{\log x}\right)\right].
\end{split}
\]
Hence, we have
\[
I(x, d)
 = e^{-\bar{d}} \frac{x}{(\log x)^{2}}\left[1 + O\left(\frac{\bar{d} \log\log x}{\log x}\right)\right],
\]
and the required estimate now follows since, if $d \sim \lambda \log x$, $\bar{d} \sim \lambda$ as $x \to \infty$. Hence, the proof of the theorem is completed.

\begin{acknowledgment}
The authors would like to express their sincere gratitude to the referee for his comments on the earlier version of this paper.
\end{acknowledgment}

\end{document}